\documentclass[a4paper,11pt]{article}
\usepackage[T1]{fontenc} 
\usepackage[cp1250]{inputenc}
\usepackage{lmodern}

\usepackage[english]{babel}
\usepackage{latexsym}
\usepackage{amsmath}
\usepackage{amsthm}
\usepackage{amssymb}
\usepackage{indentfirst}
\usepackage{amsfonts}
\usepackage{stackrel}
\usepackage{dsfont}
\usepackage{tikz}
\usepackage[mathscr]{eucal}

\newtheorem{thm}{Theorem}
\newtheorem{lemma}[thm]{Lemma}
\newtheorem{cor}[thm]{Corollary}

\newtheorem{ex}{Example}

\newtheorem{defin}{Definition}
 
\newenvironment{pro}{\begin{flushleft} \textbf{Proof}\\* \end{flushleft}}{\hfill\(\blacksquare\) \\ }

\newcommand{\reals}{\mathbb{R}}
\newcommand{\naturals}{\mathbb{N}}

\newcommand{\complex}{\mathbb{C}}

\newcommand{\eps}{\varepsilon}
\newcommand{\U}{\mathcal U}

\newcommand{\conv}{\operatorname{conv}}

\newcommand{\Ffamily}{\mathcal{F}}
\newcommand{\Gfamily}{\mathcal{G}}
\newcommand{\walt}{\operatorname{Wall}(T)}

\newcommand{\Powerset}{\mathcal{P}}

\newcommand{\Addresses}{{
  \bigskip
  \footnotesize
  Turobo\'s F., \textsc{\L\'od\'z University of Technology, \ Institute of Mathematics, \ W\'ol\-cza\'n\-ska 215, \
  	90-924 \ \L\'od\'z, \ Poland}\par\nopagebreak
  \textit{E-mail address} : \texttt{filip.turobos@gmail.com}
}}

\newcommand{\inv}{^{-1}}
\newcommand{\wall}{\operatorname{Wall}}
\newcommand{\zero}{\mathcal Z (T)}
\newcommand{\cozero}{\mathcal Z^c (T)}
\begin{document}
\title{A measure of non-compactness on $T_{3\frac 1 2 }$ spaces based on Arzel\`a-Ascoli type theorem}
\author{Filip Turobo\'s}
\maketitle

\begin{abstract}
The purpose of this paper is to generalize the measure of non-compactness for the space of continuous functions over the $T_{3 \frac{1}{2}}$ space. 
Motivated by the generalized Arzel\`a-Ascoli theorem for Tichonoff space $T$ via Wallman compactifiaction $\operatorname{Wall}(T)$, 
we constuct a measure of non-compactness for the space $C^b(T)$. We also study some of the properties of this object and give another version of Darbo-type theorem, suitable for this particular case.
\end{abstract}


\section{Introduction}

Measures of non-compactness have some interesting applications in the fixed point theory. Consequently, they play an important role in many branches of nonlinear analysis, optimization, integral and differential equations. Their usefulness stems from the fact that many problems in the mentioned fields can be reformulated in terms of a fixed point problem.

A measure of non-compactness is basically a function, defined on set of all bounded subsets of a metric space, which determines how far a set strays from being relatively compact. It takes the value of zero at each relatively compact subset of the underlying metric space.

First notes about the measures of non-compactness are associated with Kuratowski \cite{Sutherland}. In his paper \cite{Kuratowski}, he introduced the following
\begin{defin}(Kuratowski measure of non-compactness)\\
	Let $(X,d)$ be a metric space. The function $\alpha :\mathcal{P} (X) \to [0,+\infty]$ ( where $\Powerset(X)$ is a power set of $X$), defined as follows:
$$\alpha(B) = \inf \bigg\{ \delta > 0 \ : \ B \subset B_1 \cup \ldots B_n, \ \text{diam}(B_k)\leq \delta, \ 1\leq k\leq n\bigg\},$$
is called a \textbf{Kuratowski measure of non-compactness}.
\end{defin}
A similar construction which is based on $\eps$-nets has been provided in \cite{malkovski}.

\begin{defin} For a given bounded subset $A$ of a metric space $X$, its \textbf{Hausdorff measure of non-compactness} is defined as 
\[	\chi(A):= \inf\{\eps>0 \; : \; \text{there exists $\eps$-net for A}. \}
\]\end{defin}
25 years later, an Italian mathematician Darbo in \cite{Darbo} (pages 84-92) deployed the first of these functionals in his work, resulting in the celebrated Darbo fixed point theorem. Many others mathematicians followed this train of thought -- which resulted in constructing different measures of non-compactness. An interesting observation -- which turned out to be the genesis of this paper -- was made by Krukowski, who stated that whenever we have some characterization of compactness, it is possible to construct a measure of non-compactness based on this characterization.

In this paper, we present a new measure of non-compactness, originating from the refurbished version of Arzel\`a-Ascoli theorem by Krukowski \cite{ArzelaAscoli}. Some necessary definitions are provided in the end of chapter 1. Then we recall the aforementioned Arzela-Ascoli theorem to construct a new measure of non-compactness and prove some of its properties. Lastly, in chapter 3, we provide a glimpse of its possible uses.
\subsection{Notation and basic definitions}

Throughout the paper, $Y$ denotes a Banach space. $B(c,r)$ and $\overline{B}(c,r)$ are understood to be open and closed balls, respectively, centered at $c$ and of radius $r>0$. 

 For any $T_{3 \frac{1}{2}}$ space $T$, by $C(T)$ we denote set of all continuous functions on $X$, while $C^b (T)\subset C(T)$ denotes the family of continuous and bounded functions. 
 Such spaces have two families of important, distinguishable sets, namely zero and cozero sets.
 \begin{defin}(Zero and cozero sets)\\
 	Let $T$ be a Tichonoff space. By $\zero$ we denote the family
 	\[
 	\zero := \{f\inv \{0 \} \; : \; f\in C^b (T) \}.
 	\]
 	We will call set $A\in \zero$ a \textbf{zero set}. Compliments of zero sets, i.e. sets $B\subset T$ for which $T\setminus B \in \zero$ will be called \textbf{cozero sets} and the family of all such sets will be denoted by $\cozero$.
 \end{defin}
 
 These two families of sets give us a makeshift of normality in a space. This is due to the fact, that each two disjoint zero sets (which of course are closed in $T$) can be separated by two disjoint cozero sets.
 The following lemma formally represents this claim, though we will omit the proof of this simple fact.
 \begin{lemma}
 	For any $ F,G \in \zero$, $F\cap G = \emptyset$, there exist $U,V\in \cozero$ for which $F\subset U$, $G\subset V$ and $U\cap V = \emptyset$.
 \end{lemma}
 
 It turns out that zero and cozero sets are the cornerstone of the Wallman compactification construction. The elements of the compactification are called $\omega$-ultrafilters and the definition of such object is given below.
 
 \begin{defin}($\omega$-ultrafilter)\\
 	A family $\U \subset \zero$ is called a \textbf{ $\omega$-ultrafilter} if it satisfies the following conditions:
 	\begin{itemize}
 		\item[($\omega 1$)] For any finite $\Ffamily \subset \U$ its intersection is not empty, i.e. $\bigcap \Ffamily \neq \emptyset$,
 		\item[($\omega 2$)] $\U$ is maximal in the sense that it is not possible to add any set to $\U$ without losing property ($\omega 1$).
 	\end{itemize}
 \end{defin}
 
 It is convenient to think of $\omega$-ultrafilters as of maximal families of closed sets which test the compactness of a space (compare \cite{Munkres} (pages 169-170)).
 This is due to the fact, that for a non-compact space there exists a family of closed sets (which can be extended to an $\omega$-ultrafilter due to Tukey's lemma, see \cite{Jech} (page 10)), which has an empty intersection. An important characterization of $\omega$-ultrafilters is given below:
 
 \begin{lemma}
 	A family $\U \subset \zero$ is $\omega$-ultrafilter if and only if the following properties hold
 	\begin{itemize}
 		\item[(U1)] $\forall_{n\in \mathbb{N}} \left(A_1,\dots, A_n \in \U \implies \bigcap_{i=1}^n \in \U\right)$,\\
 		\item[(U2)] $ \forall_{A\in \zero} \left(\forall_{B\in \U} A\cap B \neq \emptyset \right) \implies A\in \U$.
 	\end{itemize}
 \end{lemma}
 
 As almost straightforward consequence of this characterization we obtain
 \begin{cor} \label{corol}
 	For any $A,B\subset \zero$ and any $\omega$-ultrafilter $\U$ 
 	\[
 	A\cup B \in \U \implies \left( A\in\U \; \lor \; B\in\U \right).
 	\]
 \end{cor}
 
 The very idea behind the Wallman compactification is to fill this \textit{gap} in a space with some missing elements, which this intersection lacks -- and the maximal family with this property itself turns out to be suitable for this purpose.
 
 We can now proceed to the definition of the Wallman compactification.
 
 \begin{defin}(Wallman compactification)\\
 	Let $T$ be a Tichonoff space. By $(\wall (T), \varphi)$ we will denote its \textbf{Wallman compactification}, where $\varphi: T\to \wall (T)$ is given by 
 	\[
 	\forall_{t\in T}\;  \varphi(t) := \mathcal P _t := \left\{ f\inv (0) \in \zero \; : \; f(t)=0 \right\},
 	\]
 	and $\wall (T)$ is the set of all $\omega$-ultrafilters on $T$.
 	The function $\varphi$ is sometimes called a \textbf{principal function}.
 \end{defin}

The fact that Wallman compactification is equivalent to the Cech-Stone compactification for $T_{3 \frac 1 2 }$ spaces is quite well-known. We provide the sketch of the proof for the sake of Readers' convenience.

To do this, we will firstly need a short description of topology on $\wall(T)$. We will do it by defining a base of open sets on $\wall(T)$. Fix an open set $U\subset T$ and define
\[
U^*:= \{ \U \in \wall(T) : T\setminus U \notin \U \}.
\]

It is relatively easy to see, that 
\begin{equation}\label{intersectionformula}
U^* \cap V^* = (U\cap V)^*
\end{equation}
for any $U,V\subset T$ which are open. Indeed, assume that $\U \in (U\cap V)^*$. Then $T\setminus (U\cap V) = (T\setminus U) \cup (T\setminus V) \notin \U$. If $T\setminus U$ was an element of $\U$, then (since $\U$ is maximal and $T\setminus U \subset T\setminus(U\cap V) $) $T\setminus (U\cap V) \in \U$, which contradicts that $\U \in U^*\cap V^*$. The same reasoning holds in the second case, therefore $U^*\cap V^* \supset (U\cap V)^*$.

Assume now that $\U$ belongs to both $U^*$ and $V^*$. Then $T\setminus U \notin \U$ and $T\setminus V \notin \U$. Almost analogously to the previous part of proof, if $(T\setminus U) \cup (T\setminus V) \in \U$, then from the properties of ultrafilters (namely Corollary \ref{corol}), either $T\setminus U$ or $T\setminus V$ belongs to $\U$, a contradiction.

Every ultrafilter belongs to some open set $U^*$, therefore by \cite{Engelking} (page 12) we obtain that 
$\{ U^* \; : \; U \subset T \text{ is open} \}$ form a base of topology in $\wall (T)$. This topology will be denoted by $\tau^*$.

The proof that the pair $(\wall(T), \tau^*)$ is a compactification of the Tichonoff space goes as follows. Firstly we prove that $\tau^*$ is $T_2$. Then we prove that $(\wall(T), \tau^*)$ is compact, which is done by contradiction. Lastly, we show that each continuous and bounded function defined on $T$ can be extended to a continuous function on whole $\wall(T)$ -- which proves the maximality of this compactification.
The mapping $\Gamma:C^b(T) \to C(\wall(T))$ which pairs an arbitrary function $f\in C^b(T)$ with its extension on $C(\wall(T))$ is called a \textbf{Gelfand transform}. Its properties are described in the following theorem, which is a composition of results obtained in \cite{ArzelaAscoli} and \cite{Kaniuth} (pages 66-73). These attributes of the Gelfand transform will be necessary in the latter part of our work.
\begin{thm}\label{gelfneim}
	The function $\Gamma: C^b (T) \to C\left(\wall (T)\right)$ is a homeomorphism. Moreover, it is a isometric $*$-isomorphism between these two spaces.
\end{thm}

\section{Measure of non-compactness on $C^b(T)$}\label{sectionquasimeasureofnon-compactness}

This part of the paper is dedicated to an introduction of a measure of non-compactness on $C^b(T)$, where $T$ is an arbitrary Tichonoff space. Therefore, unless explicitly stated otherwise, we will assume that $T$ is $T_{3 \frac 1 2 }$ space from now on.

The following definition introduces the concept of \textit{measure of non-compactness}. The list of axioms \textbf{(MN1)} - \textbf{(MN6)} is based on the classical collections that can be found in \cite{BanasGoebel} (page 11), \cite{Benavides} (pages 18-19) or \cite{BanasMursaleen} (page 170). 

\begin{defin}(Measure of non-compactness)\\
A function $\Omega : \Powerset(Y) \rightarrow \reals_+$ is called a \textbf{measure of non-compactness} in $Y$ if
\begin{description}
	\item[\hspace{0.4cm} (MN1)] for every $A \in \Powerset(Y), \; \Omega(A) = 0$ if and only if $A$ is relatively compact,	
	\item[\hspace{0.4cm} (MN2)] for every $A,B \in \Powerset(Y)$ such that $A \subset B$ we have $\Omega(A) \leq \Omega(B)$,
	\item[\hspace{0.4cm} (MN3)] for every $A \in \Powerset(Y)$ we have $\Omega(\overline{A}) = \Omega(A)$,
	\item[\hspace{0.4cm} (MN4)] for every $A \in \Powerset(Y)$ and $\lambda \in \complex$ we have $\Omega(\lambda A) = |\lambda|\Omega(A)$,
	\item[\hspace{0.4cm} (MN5)] for every finite $A$ and $B \in \Powerset(Y)$ such that $A \cap B = \emptyset$ we have \mbox{$\Omega(A\cup B) = \Omega(B)$},
	\item[\hspace{0.4cm} (MN6)] for every $A \in \Powerset(Y)$ we have $\Omega(\conv (A)) = \Omega(A)$.
\end{description}
\label{quasimeasureofnon-compactness}
\end{defin}

Keeping in mind that any characterisation of relatively compact sets allows us to construct an appropriate measure of non-compactness, we start with the following theorem from \cite{ArzelaAscoli}.

\begin{thm}\textbf{(Arzela-Ascoli theorem via the Wallman compactification)} \; 
The family $\Ffamily \subset C^b(T)$ is relatively compact iff
\begin{description}
	\item[\hspace{0.4cm} (AA1)] $\Ffamily$ is pointwise bounded at every point of $T$, i.e.
		$\Ffamily_t := \{f(t) \; : \; t\in T \}$ is bounded for each point $t\in T$;
	\item[\hspace{0.4cm} (AA2)] $\Ffamily$ is $\omega$-equicontinuous i.e.
	$$
	\forall_{\substack{ \eps > 0 \\ \U\in\walt}} \; \exists_{\substack{V\in \tau \\ \U\in V^*}} \; \forall_{\substack{f\in \Ffamily \\  t\in V}} \; \left|  f(t) - \hat{ f} (\U)  \right|  < \eps,
	$$
	where $\hat{f}:= \Gamma (f)$.
\end{description}
\label{AAW}
\end{thm}

To constuct the measure of non-compactness according to the theorem above, we will need to measure both pointwise boundedness and $\omega$-equicontinuity. We shall first take care of pointwise boundedness (or, since we deal with $\mathbb{C}$, relative compactness of images of each evaluation functional). It is noteworthy, that we can employ Kuratowski or Hausdorff measure of non-compactness (see \cite[Definition 2.1.]{Benavides}) in the definition below.

\begin{defin}\label{etadefin}
	Let $\xi$ be any measure of non-compactness on $\reals$. For any $\Ffamily\subset C^b (T)$ define 
	\[
	\eta (\Ffamily):= \sup_{t\in T} \xi \left( \left\{ f(t) \; : \; f\in \Ffamily \right\} \right).
	\]
	We will call functional $\eta$ a \textbf{measure of pointwise boundedness} on $C^b(T)$.
\end{defin}

In case $X = \reals$, such a functional has been employed in the study of differential equations, for example in \cite{PrzeradzkiAPM}. The functional $\eta$ measures how much a set deviates from being pointwise relatively compact, which justifies the name we used in the definition \ref{etadefin}.  The next two lemmas put this claim in a formal mathematical setting.
\begin{lemma}\label{etaprop}
	Functional $\eta$ from Definition \ref{etadefin} satisfies \textbf{(MN2)}--\textbf{(MN6)}.
\end{lemma}
\begin{pro}
Ad \textbf{(MN2)}. If $\Ffamily \subset \Gfamily \subset C^b (T)$, then -- since $\xi$ is monotone -- we obtain:
\[
\eta(\Ffamily)= \sup_{t\in T} \xi \left( \left\{ f(t)\; : \; f\in \Ffamily \right\} \right) \leqslant
\sup_{t\in T} \xi \left( \left\{ f(t)\; : \; f\in \Gfamily \right\} \right) = \eta(\Gfamily).
\]

\noindent
Ad \textbf{(MN3)}. Due to \textbf{(MN2)} we simply need to show that $\eta( \overline{\Ffamily}) \leqslant \eta (\Ffamily)$. 
Fix a point $t_0\in T$. If $f\in \overline{\Ffamily}$, then there exists a sequence $(f_n)\subset \Ffamily$ such that $(f_n)\to f$. Thus $(f_n(t_0)) \to f(t_0)$. Therefore
\[
\{f(t_0) \; : \; f\in \overline \Ffamily \} \subset \overline{ \{ f(t_0) \; : \; f\in \Ffamily \} },
\]
which, after taking closure of both these sets yields
\[
\overline{ \{f(t_0) \; : \; f\in \overline \Ffamily \} } \subset \overline{ \{ f(t_0) \; : \; f\in \Ffamily \} }.
\]
Altogether this gives us
\[
\xi \left( \left\{ f(t_0) \; : \; f\in \overline{\Ffamily} \right\} \right) \leqslant \xi \left( \overline{ \left\{ f(t_0) \; : \; f\in \Ffamily \right\}  } \right) = \xi \left( { \left\{ f(t_0) \; : \; f\in \Ffamily \right\}  } \right)
,\]
where the last equality follows from \textbf{(MN3)} property of $\xi$. Since $t_0$ was chosen arbitrarily, then taking supremum over all points of $T$ yields \textbf{(MN3)} property for $\eta$.

\noindent
Ad \textbf{(MN4)}. Let $\lambda \in \mathbb C$. We have that
\begin{eqnarray*}
\eta (\lambda\Ffamily) &=& \sup_{t\in T} \xi \left( \left\{ \lambda f(t) \; : \; f\in \Ffamily \right\} \right) \\
&=& \sup_{t\in T} \xi \left( \lambda \left\{ f(t) \; : \; f\in \Ffamily \right\} \right) \\
&=&\sup_{t\in T} |\lambda| \xi \left( \left\{ f(t) \; : \; f\in \Ffamily \right\} \right) \\
&=& \lambda \eta (\Ffamily).
\end{eqnarray*}

\noindent
Ad \textbf{(MN5)}. Let $\Gfamily$ be a finite family of continuous and bounded functions on $T$. Therefore, from \textbf{(MN5)} property of $\xi$ we obtain
\begin{eqnarray*}
\eta(\Ffamily \cup \Gfamily) &=& \sup_{t\in T} \xi \left( \{f(t) \; : \; f\in \Ffamily \} \cup \{g(t) \; : \; g\in \Gfamily \}\right) \\
&=& \sup_{t\in T} \xi \left(\{f(t) \; : \; f\in \Ffamily \} \right)\\
&=& \eta (\Ffamily).
\end{eqnarray*}

Ad \textbf{(MN6)} Since $\beta$ is a measure of non-compactness on $\complex$, it is invariant under passing to the convex hull.
 Notice that for a fixed $t\in T$ and $\Ffamily \subset C^b(T)$  the following sets are equal:
 \[
 \left\{ f(t) \; : \; f\in\conv(\Ffamily) \right\} = \conv \left( \left\{ f(t) \; : \; f\in \Ffamily \right\} \right).
 \]
 Combining these two facts together yields
 \begin{eqnarray*}
 	\eta\big(\conv(\Ffamily)\big) &=& \sup_{t\in T} \beta \left( \left\{ f(t) \; : \; f \in \conv(\Ffamily) \right\} \right)\\
 	\;&=& \sup_{t\in T} \beta \left( \conv\left( \left\{ f(t) \; : \; f \in \Ffamily \right\} \right) \right)\\
 	\;&=& \sup_{t\in T} \beta \left(\left\{ f(t) \; : \; f \in \Ffamily \right\} \right) \\
 	\;&=& \eta(\Ffamily).
 \end{eqnarray*}
 Thus, $\eta$ has the desired invariance property.
\end{pro}
The next lemma justifies clearly, why we refer to $\eta$ as to a measure of pointwise boundedness (compactness).
\begin{lemma}\label{etaequi}
	For $\Ffamily \subset C^b(T)$ the following equivalence holds:
	
	{\begin{center}
	$\Ffamily$ is pointwise bounded $\iff$ $\eta(\Ffamily) = 0$.
\end{center}
	}
\end{lemma}
\begin{pro}
	
	"$\implies$" If $\Ffamily$ is pointwise bounded, then for any $t\in T$ the set $\{f(t) \; : \; f\in \Ffamily\}$ is bounded. Hence
$\xi \left( \{f(t) \; : \; f\in \Ffamily\} \right) = 0$. Taking supremum over all $t\in T$ yields the desired implication.

	"$\impliedby$" If $\eta(\Ffamily) = 0$, then for any $t\in T$ the set $\{f(t) \; : \; f\in \Ffamily\}$ is relatively compact, hence bounded. Therefore, $\Ffamily$ is pointwise bounded.
\end{pro}

Apart from measuring pointwise relative compactness, we would like to measure the violation of $\omega$-equicontinuity. The following, stepwise construction will provide us with a proper tool for that purpose.
\begin{defin}\label{Omegadefinition}
	For every $\Ffamily \subset C^b(T)$, $\U \in \wall (T)$ and open $V\subset T$, for which $\U \in V^*$ we define
	\[
	\omega^\U (\Ffamily, V):= \sup_{f\in \Ffamily} \sup_{t\in V} \left| \hat{f}(\U) - f(t) \right|;
	\]
	\[
	\omega^\U (\Ffamily) := \inf_{V^* \ni \U} \omega^\U (\Ffamily, V);
	\]
	\[
	\omega(\Ffamily):= \sup_{\U\in \wall X} \omega^\U (\Ffamily).
	\]
\end{defin}

The functional $\omega$ defined above will measure the violation of $\omega$-equicontinuity. To handle some difficulties with proving one of its properties, we will need an additional
\begin{lemma}\label{restriction}(Restriction Lemma)
	Let $\Ffamily\subset C^b(T)$, $\U \in \walt$ and $V\subset T$ be a fixed open subset, for which $\U \in V^*$. The following equality holds:
	\[
	\omega^\U (\Ffamily) = \inf_{\stackrel{\tau \ni W\subset V}{W^* \ni \U}} \omega^\U (\Ffamily, V).
	\]
\end{lemma}
\begin{pro}
	Let $S:=\{W\in \tau \; : \; \U\in W^* \subset V^* \}$ and $S':= \{W\in \tau : \U \in W^* \not\subset V^*\}$. Note that
\[\inf_{\U \in W\in \tau} \omega^\U (F,W) = \min \left\{ \; \inf_{W\in S} \omega^\U (F,W),\; \inf_{W\in S'} \omega^\U (F,W) \;\right\} .\]
It is enough to prove that
\[
\inf_{W\in S} \omega^\U (F,V) \leqslant \inf_{W\in S'} \omega^\U (F,V).
\]
Take any set $W\in S'$. Since $\U \in W^*$, then $W\cap V$ is non-empty, due to \eqref{intersectionformula} and the simple fact, that $\emptyset^* = \emptyset$. Therefore $W\cap V \in S$. Also notice that $\omega^\U (F,W\cap V) \leqslant \omega^\U (F,W)$, since we take supremum over the smaller (in the sense of inclusion) set. 

Therefore, for each $W\in S'$ we have a corresponding set of form $(V\cap W)\in S$, on which the functional $\omega^\U$ attains a smaller value. Thus, we have proven that
\[
\inf_{W\in S} \omega^\U (F,V) \leqslant \inf_{W\in S'} \omega^\U (F,V).
\]
\end{pro}

Now, we will show that our functional $\omega$ has all the properties \textbf{(MN2)}-\textbf{(MN6)}, which is the statement of the following

\begin{lemma}\label{omegaprop}
The functional $\omega$ satisfies the properties \textbf{(MN2)}--\textbf{(MN6)}.
\end{lemma}
\begin{pro}
	Ad \textbf{(MN2)}. Let $\Ffamily \subset \Gfamily \subset C^b (T)$ and fix $\U \in \wall (T)$. Let $V$ be any open subset of $T$ for which $\U\in V^*$. Therefore 
	\[
	\omega^\U (\Ffamily, V) = \sup_{f\in \Ffamily} \sup_{t\in V} \left| \hat{f}(\U) - f(t) \right| \leqslant \sup_{f\in \Gfamily} \sup_{t\in V} \left| \hat{f}(\U) - f(t) \right| = \omega^\U (\Gfamily, V).
	\]
	Since $V$ was arbitrary, then $\omega^\U(F) \leqslant \omega^\U (G)$. Due to $\U$ being arbitrary, we obtain $\omega(\Ffamily) \leqslant \omega(\Gfamily)$, thus proving the first property.
	
	\noindent Ad \textbf{(MN3)}. Fix $\eps>0$. Once again consider a point $\U \in \wall (T)$ and its open neighbourhood $V^*$. From the definition of supremum, there exists a function $f_1 \in \overline{\Ffamily}$ such that
	\[
	\omega^\U (\overline{\Ffamily}, V) = \sup_{f\in \overline{\Ffamily}} \sup_{t\in V}\left|\hat{f}(\U) - f(t) \right| \leqslant \sup_{t\in V} \left| \hat{f_1}(\U) - f_1 (t)  \right| + \frac \eps 3. 
	\]
	Since $f_1 \in \overline \Ffamily$, every neighbourhood of $f_1$ intersects $\Ffamily$ in a non-empty way. The mapping $\Gamma$ is an isometry -- hence for any positive constant, then in particular for $\eps >0$, the following equivalence holds
	\[
 \left(	\sup_{t\in T} \left| f_1(t) - f_2 (t) \right| < \frac{\eps}{3}  \right) \iff \sup_{\mathcal V\in \wall (T)} \left|\hat{f_1}(\mathcal V) - \hat{f_2}(\mathcal V)\right| < \frac \eps 3
	\]
	for any $f_2\in C^b(T)$. The intersection of $B\left(f_1,  \frac \eps 3 \right)$ and $\Ffamily$ is non-empty, so let $f_2\in B(f_1,  \frac \eps 3 )\cap \Ffamily$. Then
	\[
\left( \forall_{t\in T} \;\; \left|f_1(t) - f_2(t) \right| <\frac \eps 3  \;  \right) \; \land \; \left(\; \left| \hat{f_1} (\U) - \hat{f_2}(\U) \right| < \frac \eps 3 \; \right).  
	\]
	In conclusion, we obtain
	\begin{eqnarray*}
	\omega^\U (\overline{\Ffamily}, V)& \leqslant &\sup_{t\in V} \left| \hat{f_1}(\U) - f_1 (t)  \right| + \frac \eps 3\\
	&\leqslant & \sup_{t\in V} \left( \left| \hat{f_1}(\U) - \hat{f_2} (\U)  \right| + \left|\hat{{f_2}}(\U) - f_2(t) \right| + \left| f_2(t) - f_1(t)\right| \right) + \frac \eps 3\\
	&\leqslant & \frac \eps 3 + \sup_{t\in V} \left|\hat{{f_2}}(\U) - f_2(t) \right| + \frac \eps 3 + \frac \eps 3 \\
	&\leqslant & \omega^{\U} (\Ffamily, V) + \eps.
	\end{eqnarray*}
	Since $\varepsilon$ was arbitrary, we obtain the desired property.
	
	\noindent Ad \textbf{(MN4)}. This property is obvious due to homogenity of $\Gamma$.
	\[
	\omega^\U(\lambda \Ffamily, V) := \sup_{f\in \Ffamily} \sup_{t\in V} \left| \widehat{\lambda f}(\U) - \lambda f(t) \right| =
	\sup_{f\in \Ffamily} \sup_{t\in V} |\lambda |\cdot \left| \hat{ f}(\U) - f(t) \right| = |\lambda | \omega^\U(\Ffamily,V).
	\]
	Passing to the infimum over all neighbourhoods of $\U$ and then taking supremum over all $\U \in \wall (T)$ yields the desired property.
	
	\noindent Ad \textbf{(MN5)}. The proof of this property can be easily done by induction, with additional help of Restriction Lemma \ref{restriction}. We will prove that extending any $\Ffamily \subset C^b(T)$ by any finite set of functions from this space will not change the value of $\omega$ on this family.
	
	Therefore, let $\Ffamily \subset C^b(T)$ and $f_0\notin \Ffamily$ be another function from $C^b(T)$. Fix $\U \in \walt$ and a positive $\eps$. We will prove that 
	\[
	\omega^\U \left(\Ffamily \cup \{f_0 \}\right) = \omega^\U(\Ffamily) + \eps.
	\]
	The singleton $\{f_0\}$ is trivially $\omega$-equicontinuous, therefore we are able to find $V\in \tau$ such that $\U\in V^*$ and 
	\begin{equation}\label{epsilonek}
	\forall_{t\in V} \; \; \left| f_0(t) - \widehat{f_0}(\U)  \right|  < \eps.
	\end{equation}
	Now, due to the Restriction Lemma, if we denote by $S:=\left\{
	W\in \tau \; : \; \U\in W^* \subset V^*
	\right\}$ then, taking into consideration \eqref{epsilonek}, we obtain
	\[\begin{array}{cclcl}
	\omega^\U \left(\Ffamily \cup \{f_0 \} \right) &=&  \inf_{W\in S}  \omega^\U 
	\left(\Ffamily\cup \{f_0\}, W\right) &=&  \inf_{W\in S} \max \{ \omega^\U 
	\left(\Ffamily, W\right), \eps \}
	\\ \; &\leqslant& \inf_{W\in S} \omega^\U \left(\Ffamily, W\right) +\eps &=& \omega^\U \left(\Ffamily \right) + \eps
	\end{array}\]
	Since $\eps>0$ was arbitrary, then $\omega^\U\left(\Ffamily \cup \{f_0\}\right) = \omega^\U\left( \Ffamily \right)$. Due to $\U\in\walt$ being arbitrary as well, we obtain that
	 $\omega \left( \Ffamily \cup \{f_0\} \right) = \omega (\Ffamily)$.
	Inductively, we can add any finite family of functions from $C^b(T)$ -- which proves the fifth desired property of $\omega$.
	
	This proves that $\omega$ is a quasimeasure of non-$\omega$-equicontinuity.
	
	Ad \textbf{(MN6)}. After brief examination of a problem, one might see a necessity for referring to the characterisation of convex hull.
	\begin{lemma}
		Let $X$ be a linear space and let $A\subset X$. The following equivalence holds
		
		{\begin{center}
		$x\in \conv(A)$ $\iff$ $x$ is a convex combination of elements of $A$.
	\end{center}}
	
	\end{lemma}  
	Taking this fact into account, if we fix $\U$ and an open set $V\subset T$ for which $\U \in V^*$ we have the following equalities:
	\begin{eqnarray*}
		\omega^\U \big(\conv(\Ffamily), V\big) &=& \sup_{g\ \in \conv (\Ffamily)} \sup_{t\in V}  \left| g(t) - \hat{g}(\U) \right|\\
		\;&=& \sup_{n\in \mathbb N} \sup_{\substack{s_1,\dots s_n >0\\ \sum_{i=1}^n s_i =1 }} \sup_{f_1,\dots, f_n \in \Ffamily}  \sup_{t\in T} \left| \sum_{i=1}^n s_i f_i(t) - \sum_{i=1}^n \widehat{s_i f_i}(\U) \right|\\
	\end{eqnarray*}
	Due to linearity of the Gelfand transform we can split the extension of convex combination of the right fragment of the difference, obtaining
	\begin{eqnarray*}
		\omega^\U \big(\conv(\Ffamily), V\big) &=&\sup_{n\in \mathbb N} \sup_{\substack{s_1,\dots s_n >0\\ \sum_{i=1}^n s_i =1 }} \sup_{f_1,\dots, f_n \in \Ffamily}  \sup_{t\in T} \left| \sum_{i=1}^n s_i f_i(t) - \sum_{i=1}^n s_i\hat{ f_i}(\U) \right|\\
		\;&\leqslant&\sup_{n\in \mathbb N} \sup_{\substack{s_1,\dots s_n >0\\ \sum_{i=1}^n s_i =1 }} \sup_{f_1,\dots, f_n \in \Ffamily}  \sup_{t\in T} \; \sum_{i=1}^n s_i \left| f_i(t) - \hat{ f_i}(\U) \right|\\	
		\;&\leqslant&\sup_{n\in \mathbb N} \sup_{\substack{s_1,\dots s_n >0\\ \sum_{i=1}^n s_i =1 }} \sup_{f_1,\dots, f_n \in \Ffamily} \; \sum_{i=1}^n s_i \cdot \left(  \sup_{t\in T}  \left| f_i(t) - \hat{ f_i}(\U) \right| \right), \\
	\end{eqnarray*}
	since we can maximize each addend independently. By maximizing this expression through each function $f_i$ separately we obtain
	\begin{eqnarray*}
		\omega^\U \big(\conv(\Ffamily), V\big) &\leqslant&\;\sup_{n\in \mathbb N} \sup_{\substack{s_1,\dots s_n >0\\ \sum_{i=1}^n s_i =1 }} \; \sum_{i=1}^n s_i \cdot \left( \sup_{f_i \in \Ffamily}   \sup_{t\in T}  \left| f_i(t) - \hat{ f_i}(\U) \right| \right) \\
		&=&\;\sup_{n\in \mathbb N} \sup_{\substack{s_1,\dots s_n >0\\ \sum_{i=1}^n s_i =1 }} \; \sum_{i=1}^n s_i \cdot \left( \sup_{f \in \Ffamily}   \sup_{t\in T}  \left| f(t) - \hat{ f}(\U) \right| \right) \\
		&=&\;  1 \cdot \left( \sup_{f_i \in \Ffamily}   \sup_{t\in T}  \left| f_i(t) - \hat{ f_i}(\U) \right| \right) \\
		\;	&=&\omega^\U (\Ffamily, V).
	\end{eqnarray*}
	Since both $V$ and $\U$ were arbitrary, we get
	\begin{equation*}
	\omega^\U \big(\conv(\Ffamily)\big) \; = \;	\inf_V \omega^\U \big(\conv(\Ffamily), V\big) \; \leqslant \; \inf_V \omega^\U \big(\Ffamily, V\big) \; = \; \omega^\U (\Ffamily).
	\end{equation*}	
	Finally, $\omega\big( \conv( \Ffamily ) \big) \leqslant \omega (\Ffamily)$. Since $\omega$ satisfies \textbf{(MN2)} and $\conv(\Ffamily) \supset \Ffamily$, then $\omega$ is invariant under passing to convex hull.
\end{pro}

The fact that $\omega$ measures how much $\Ffamily \subset C^b(T)$ violates $\omega$-equicontinuity is expressed by the lemma below:
\begin{lemma}\label{omegaequi}
	For any family $\Ffamily \subset C^b(T)$, the following equivalence holds:
	\begin{center}
		$\Ffamily$ is $\omega$-equicontinuous $\iff$ $\omega(\Ffamily) = 0$.
	\end{center}
\end{lemma}
\begin{pro}
	"$\implies$" Assume that $\Ffamily$ is $\omega$-equicontinuous. Take any $\U \in \wall (T)$. Then, for every $\eps_n :=\frac 1 n$ there exists $V_n \subset T$ -- open, such that for each $f\in \Ffamily$, $t\in V_n$ the inequality $|\hat{f}(\U) - f(t)|<\frac 1 n$. Hence
	\[
	\omega^\U (\Ffamily, V_n ) \leqslant \frac 1 n,
	\]
	and since $\U \in V_n^*$ for all $n\in \mathbb N$, $\omega^\U(\Ffamily) = 0$. Since $\U$ was chosen arbitrarily, $\omega(\Ffamily) =0$.
	
	"$\impliedby$" Assume that $\omega (\Ffamily) = 0$. Therefore, for all $\U \in \wall (T)$ we have $\omega^\U (\Ffamily) = 0$. For the sake of contradiction, suppose that $\Ffamily$ fails to be $\omega$-equicontinuous. Then, there exists $\eps>0$ and $\U_0 \in \wall (T)$ such that for any $V^*$-neighbourhood of $\U_0$ there exists $f_V \in \Ffamily$ and $t_V\in V$ for which
	\[	\left|  \hat{f_V} (\U_0) - f(t_V) \right| > \eps.	\]
	As a result
	\[
	\left(	\forall_{V^*\ni \U_0}\; \omega^{\U_0} (\Ffamily, V) \right) > \eps \implies \omega^{\U_0}(\Ffamily) \geqslant \eps.
	\]
	This implies that $\eps \leqslant \omega^{\U_0} (\Ffamily) \leqslant \omega (\Ffamily) = 0$, a contradiction. Hence $\Ffamily$ is $\omega$-equicontinuous.
\end{pro}

Having all the necessary tools, we can present the main result of our divagations until now.

\begin{thm}\label{Omegathm}
	The functional $\Omega : C^b(T) \to [0,+\infty]$ given by the formula
	\[
	\Omega (\Ffamily):= \omega(\Ffamily) + \eta(\Ffamily)
	\]
	satisfies properties \textbf{(MN1)}-\textbf{(MN6)}.
	Therefore $\Omega$ is a measure of non-compactness on $C^b(T)$.
\end{thm}
\begin{pro}
	Thanks to Lemmas \ref{etaprop} and \ref{omegaprop}, we obtain the properties \textbf{(MN2)}--\textbf{(MN6)} immediately. 
	If $\Ffamily$ is relatively compact, then due to Theorem \ref{AAW} it is both pointwise bounded and $\omega$-equicontinuous. Consequently, due to Lemmas \ref{etaequi} and \ref{omegaequi}
	\[
		\Omega(\Ffamily) = \eta(\Ffamily) + \omega(\Ffamily) = 0.
	\]
	Conversely, if $\Omega(\Ffamily) = 0$, then both $\eta(\Ffamily)$ and $\omega(\Ffamily)$ equal to $0$. This means $\Ffamily$ is both $\omega$-equicontinuous and pointwise bounded, thus relatively compact.
\end{pro}
\begin{ex}
	Consider a compact unit interval $[0,1]$ with the topology induced by Euclidean metric and a family of continuous and bounded functions $\Ffamily:=\left\{ x^n \; : \; n\in \mathbb{N} \right\}$.\\
	One can easily see that the Cech-Stone compactification of this interval is the interval $[0,1]$ itself, since compactification of a compact space does not alter its structure. So one can identify $U\in \walt$ simply as elements of unit interval. Since $\Ffamily$ is bounded, then $\eta(\Ffamily)=0$. But this family fails to be compact, therefore it has to be non-$\omega$-continuous to some extent. Let us measure to what extent this property fails to be satisfied. Fix some $t_0 \in [0,1]$. Due to Lemma \ref{restriction}, we can restrict only to the neighbourhoods of the form $(t_0-\eps, t_0+\eps)$. In case when $t_0 = 1$ or $t_0=0$ we can use open sets of the form  $(1-\eps,1]$ and $[0,\eps)$ respectively.
		For $t_0 = 0$ and its $\eps$-neighbourhood obviously $\omega^{t_0} (\Ffamily, V) = \eps$, since $\eps<1$. And for $1$ we get
		\[
		\omega^{1} (\Ffamily, V) := \sup_{n\in\mathbb{N}} \left|(1-\eps)^n - 1^n\right| = 1,
		\] since $\lim\limits_{n\to\infty} \left|(1-\eps)^n - 1^n\right| = 1$.
	For some $t_0\in (0,1)$ and $\eps>0$ such that $V:=(t_0-\eps, t_0+\eps) \subset [0,1]$ we have
	\[
	\omega^{t_0} (\Ffamily, V) := \sup_{n\in\mathbb{N}} \sup_{t\in V} \left|x^n - (t_0)^n\right|
	\]
	Since the derivation of $(x^n - t_0^n)$ yields $x^{n-1}$, we can easily see that this supremum in this case can be calculated as the value of this expression at one of the points $(t_0-\eps)$ or $(t_0+\eps)$. Of course the latter value yields a greater difference for each $n\in \mathbb{N}$, thus we stick to the latter one.
	Thus $\omega^{0} (\Ffamily, V) = \sup_{n\in\mathbb{N}} \left((t_0 +\eps)^n - (t_0)^n \right) $. This value does not exceed $1$ though.
	If we tend with $\eps$ to $0$, i.e. we take the infimum over all neighbourhoods of $t_0\in [0,1]$, we obtain that $\omega^{t_0}(\Ffamily)\leqslant 1$ for all $t_0$ and it attains value of $1$ for $t_0 = 1$. Consequently
	$\omega(\Ffamily) =1$.\\
	Let us try to evaluate the value of the Hausdorff measure of non-compactness on this family. Obviously, this measure does not exceed $1$, since a closed ball centered at zero function contains the whole family $\Ffamily$. Notice, however, that for any $n\in\mathbb{N}$
	\[
	\sup_{n\in \mathbb{N}} \|x^n - \frac{1}{2}\| = \sup_{n\in \mathbb{N}} \sup_{x\in[0,1]} |x^n - \frac{1}{2}| = \frac{1}{2}. 
	\]
	This means that the Hausdorff measure of $\Ffamily$ cannot exceed $\frac{1}{2}$.
	
	We will show, that $\alpha (\Ffamily) \geqslant 1$. On the contrary, suppose that $\Ffamily$ can be covered by a finite family of sets $\Gfamily$, whose diameters are no greater than $1-\eps$.
	Notice that for $n<k$
	\[
	1\leqslant	\|x^n - x^k\| = \sup_{x\in[0,1]} (x^n - x^k) \geqslant \left(\frac{n}{k} \right)^{\frac{n}{k-n}} - \left(\frac{n}{k} \right)^{\frac{k}{k-n}} \stackrel{k\to\infty}{\longrightarrow} 1 - 0 = 1.
	\] 
	If $\Gfamily$ covers entire $\Ffamily$, then the identity mapping $x^1$ belongs to some $A_1 \in \Ffamily$. Since $\|x^n - x\| \stackrel{n\to\infty}{\longrightarrow} 1$, there exists $n_1\in \naturals$ such that $\|x^{n} - x\| > 1 - \eps$ for all $n\geqslant n_1$. Therefore $x^{n_1} \in A_2$ for some different $A_2\in \Gfamily$. Again, there exists $ n_1<{n_2}\in \naturals$ such that $x^n \notin A_1\cup A_2$ for all $n\geqslant n_2$. Consequently, if $\bigcup_{i=1}^k A_i$ contains all functions $x^n$ up to some $n_{k-1}$, there exists $n_k$ such that $x^n\notin \bigcup_{i=1}^k A_i$ for $n\geqslant n_k$. Thus $\Gfamily$ cannot be finite. This proves, that $\alpha(\Ffamily) = 1$. This reasoning also shows that $\chi(\Ffamily) = \frac{1}{2}$.
\end{ex} 
%

Due to the fact, that $\Gamma$ is an isometry, one may think that neither Kuratowski's nor Hausdorff's measure can exceed $\Omega$. The example provided proves that the opposite relation cannot hold in case of $\chi$. Thus we pose the following

\textbf{Question 1:}
Prove or disprove if there exists any form of inequality between $\Omega$ and some other measures of non-compactness, including Kuratowski's and Hausdorff's measures.

\section{Applications -- Darbo fixed point theorem.}
\label{applications}

The theorem we present below is stated in the classic terms of $Y$ being any Banach space and $\Omega$ being any quasimeasure of non-compactness. However, space $Y = C^b(X)$ is, in fact, a Banach space and we can implement in this theorem the measure of non-compactness constructed in Theorem \ref{Omegathm}. The proof is inspired by \cite{Lesniak} and \cite{Aghajani}, especially by the use of consecutive iterations of the function $\psi$. Nevertheless we dish it up for the Readers' convenience.
  
\begin{thm}
Let $\Omega : \Powerset(Y) \rightarrow \reals_+$ be a quasimeasure of non-compactness, $C \subset Y$ a non-empty, convex, bounded and closed subset and $\psi : \reals_+ \rightarrow \reals_+$ a nondecreasing function such that $\lim_{n\rightarrow \infty} \psi^{(n)}(t) = 0$ for every $t \geq 0$. If $\Phi : C \rightarrow C$ is a continuous function such that  

\begin{description}
	\item[\hspace{0.4cm} (B)] for all $A \subset C$ we have 
	$$\Omega(\Phi(A)) \leq \psi(\Omega(A)),$$
\end{description}

\noindent
and $\Omega(C)<\infty$ then $\Phi$ has at least one fixed point in $C$. 
\label{Darbotheorem} 
\end{thm}
\begin{pro}

Denote $C_1 := C$ and define $C_{n+1} := \overline{\conv \Phi(C_n)}$. From this definition, the closedness of $C_n$ for every $n\in\naturals$ is obvious. By induction we will prove that $(C_n)_{n\in\naturals}$ is a descending sequence. Indeed, since $\Phi(C_1)\subset C_1$, then $$C_2 :=\overline{\conv \Phi(C_1) } \subset C_1, $$ due to $C$ being a closed and convex set. Notice, that $\Phi_{|C_2} :C_2 \to C_2$, because $\Phi(C_2) \subset \Phi(C_1) \subset C_2$. Thus, using the same argument as previously, the inductive step can be proved.

If we define $D$ as $D := \bigcap_{n\in\naturals} C_n$, we see that it is an intersection of closed and convex sets, thus it is a closed and convex set itself.

By \textbf{(B)}, we have the following bound on $\Omega(C_n)$
\begin{gather}
\Omega(C_{n+1}) \stackrel{\textbf{(MN3)}}{=} \Omega(\Phi(C_n)) \leq \psi(\Omega(C_n)) \leq \psi^{(2)}(\Omega(C_{n-1})) \leq \ldots \leq \psi^{(n)}(\Omega(C)).
\label{maindarboinequality}
\end{gather}
This (paired with the assumption that $\Omega(C)<\infty$) implies that the sequence $\Omega(C_n)$ tends to $0$ for $n\rightarrow \infty$. As a consequence $\Omega(D) = 0$, thus $D$ is relatively compact, though we do not know at this point whether it contains any elements or not.

In the next part of the proof, we would like to show that $D$ is non-empty. Since each $C_n$ is non-empty, we can consider a sequence $(x_n)_{n\in\naturals}\subset C$ such that $x_n\in C_n$, then for a fixed $k \in \naturals$ we have   
$$\Omega((x_n)_{n\in\naturals}) = \Omega\bigg((x_n)_{n=1}^k \cup (x_{k+n})_{n\in\naturals}\bigg) \stackrel{\textbf{(MN5)}}{\leq} \Omega(C_{k+1}).$$

\noindent
As $k\rightarrow\infty$, we conclude that $\Omega((x_n)_{n\in\naturals}) = 0$, i.e. $\{x_n \; : \; n\in\naturals\}$ is relatively compact, thanks to the property \textbf{(MN1)}. Therefore, the existence of a convergent subsequence of $(x_n)_{n\in\naturals}$ is established. Let us denote its limit by $x\in Y$.

Now, each of the sets $C_n$ was closed, and due to the fact, that $(C_n)_{n\in \naturals}$ is descending, for a fixed positive $k\in \mathbb{N}$, the sequence $(x_{k+n})_{n\in\naturals}$ is a convergent sequence of elements from $C_k$. Since $k$ was arbitrary, $x\in C_k$ for each $k\in \naturals$, hence $x\in D$. Now, since $D$ is non-empty, relatively compact and closed, it is a non-empty and compact set, which happens to be convex as well. 

If we restrict $\Phi$ to set $D$, we obtain a continuous map, which is compact as well.

Thus, applying Schauder fixed point theorem to $\Phi :D \to D$, we can  establish the existence of the fixed point.
\end{pro}

\subsection*{Final remarks}
	From the lecture of \cite{KrukowskiPrzeradzki} it seems to be possible to create a different measure of non-compactness over the space of functions which map $\sigma$-locally compact Hausdorff space $X$ to some metric space $Y$. The result obtained by the authors of this publication is also a result in the spirit of Arzel\'a-Ascoli theorem, which seems to be promising for further research. Another direction, provided by one of the authors of the previously mentioned paper in \cite{Arzela-Ascolitheoreminuniformspaces} also contains results which are quite similar to the ones used in the construction of the functional $\Omega$. It allows us to hope that other measures of non-compactness, based on the other variants of Arzel\'a-Ascoli theorem, can possess some interesting properties. 
	
\subsection*{Acknowledgments}
	The author would like to express his utmost gratitude to his friend and auxiliary supervisor, Mateusz Krukowski, whose constructive criticism has been an invaluable help throughout the work on this paper. 
	
	Moreover, the author wishes to acknowledge the support of his supervisor Jacek Jachymski and colleague Jarosław Swaczyna. Their comments and discussion over the obtained results were truly fruitful. Also, mathematical and linguistic remarks of Katarzyna Chrząszcz led to an overall increase of the quality of this paper.

\Addresses
\end{document}